\documentclass[12pt]{article}
\newtheorem{theorem}{Theorem}
\title{The variance of the Stirling cycle numbers}
\author{Herbert S. Wilf\\
University of Pennsylvania\\
Philadelphia, PA 19104-6395}

\begin{document}
\maketitle
\begin{abstract}
We show that the probability that two permutations of $n$ letters have the same number of cycles is
\[\sim \frac{1}{2\sqrt{\pi\log{n}}}.\]
\end{abstract}

\vspace{.5in}

Our purpose here is to prove the following
\begin{theorem}
Let two permutations of $n$ letters be chosen independently uniformly at random. The probability that they have the same number of cycles is
\[\sim \frac{1}{2\sqrt{\pi\log{n}}}\qquad (n\to\infty).\]
\end{theorem}
This question was raised by Mikl\'os B\'ona.

The Stirling cycle numbers ${n\brack k}$ are defined by
\begin{equation}
\label{eq:gf}
\prod_{j=0}^{n-1}(x+j)=\sum_{k=1}^n{n\brack k}x^k,
\end{equation}
and, as is well known, ${n\brack k}$ is the number of permutations of $n$ letters that have exactly $k$ cycles. It follows that
\[\sum_k{n\brack k}=n!,\]
of course, but what can be said about $f(n)=\sum_k{n\brack k}^2$?

Since for any polynomial $g(z)=\sum_{k=0}^{n-1}a_kz^k$ with real or complex coefficients we have
\[\sum_{k=0}^{n-1}|a_k|^2=\frac{1}{2\pi}\int_0^{2\pi}|g(e^{i\theta})|^2d\theta ,\]
we have in particular that
\begin{eqnarray*}
f(n)&=&\sum_k{n\brack k}^2\\
&=&\frac{1}{2\pi}\int_0^{2\pi}\left|\prod_{j=0}^{n-1}(e^{i\theta}+j)\right|^2d\theta\\
&=&\frac{(n-1)!^2}{2\pi}\int_0^{2\pi}\left|\prod_{j=0}^{n-1}\left(1+\frac{e^{i\theta}}{j}\right)\right|^2d\theta\\
&=&\frac{(n-1)!^2}{2\pi}\int_0^{2\pi}e^{2H_{n-1}\cos{\theta}}\left|\prod_{j=0}^{n-1}\left(1+\frac{e^{i\theta}}{j}\right)e^{-e^{i\theta}/j}\right|^2d\theta ,
\end{eqnarray*}
where $H_{n-1}=\sum_{m=1}^{n-1}1/m$ is the $(n-1)$st harmonic number.

Now the familiar Gamma function of analysis is given by
\[\frac{e^{-\gamma z}}{z\Gamma(z)}=\prod_{r=1}^{\infty}\left(1+\frac{z}{r}\right)e^{-z/r},\]
and consequently
\begin{eqnarray}
f(n)&=&\frac{(n-1)!^2}{2\pi}\int_0^{2\pi}e^{2H_{n-1}\cos{\theta}}\left\{\left|\frac{e^{-\gamma e^{i\theta}}}{e^{i\theta}\Gamma(e^{i\theta})}\right|^2+o(1)\right\}d\theta \nonumber\\
&=&(1+o(1))\frac{(n-1)!^2}{2\pi}\int_0^{2\pi}\frac{e^{2\cos{\theta}\log{n}}}{|\Gamma(e^{i\theta})|^2}d\theta\nonumber\\
&=&(1+o(1))\frac{n!^2}{2\pi}\int_0^{2\pi}\frac{e^{2(\cos{\theta}-1)\log{n}}}{|\Gamma(e^{i\theta})|^2}d\theta\label{eq:intg}
\end{eqnarray}

The method of Laplace for integrals (e.g., \cite{db}) now shows that the behavior of our integral
\[I(n)=\int_0^{2\pi}\frac{e^{2(\cos{\theta}-1)\log{n}}}{|\Gamma(e^{i\theta})|^2}d\theta\]
for large $n$ is given by
\[I(n)\sim \sqrt{\frac{\pi}{\log{n}}}.\]
The final result is that
\begin{equation}
\label{eq:fnl}
f(n)\sim \frac{n!^2}{2\sqrt{\pi\log{n}}}.
\end{equation}
which completes the proof of the theorem.

\end{document}